%% file: main.tex
\documentclass{amsart}
\usepackage[utf8]{inputenc}

\usepackage{multirow}
\usepackage{verbatim}
\usepackage{color,graphicx}
\usepackage{amsmath,amsfonts,amsthm,amscd,upref,amstext}

\usepackage{import}
\usepackage{xifthen}
\usepackage{pdfpages}
\usepackage{transparent}

\newcommand{\D}{\mathbb{D}}

\newtheorem{theorem}{Theorem}[section]
\newtheorem{lemma}[theorem]{Lemma}

\newtheorem{corollary}[theorem]{Corollary}
\newtheorem{proposition}[theorem]{Proposition}

\newtheorem{remark}[theorem]{Remark}

\begin{document}

\title[Indecomposable continua \& Julia sets of polynomial-like mappings]{Indecomposable continua and the Julia sets of polynomial-like mappings}
\author{Elena Gomes}

\address{ Centro de Matemática, Facultad de Ciencias,  Universidad de la República, Uruguay.}
\email{egomes@cmat.edu.uy}

\maketitle

\begin{abstract}
    Let $f$ be a polynomial-like mapping of the sphere of degree $d \geq 2$. We show that the Julia set $J(f)$ of $f$ cannot be the union of a finite number of indecomposable proper subcontinua. As a corollary, we prove that $J(f)$ is an indecomposable continuum if and only if there exists a prime end of some complementary region of $J(f)$ whose impression is the entire $J(f)$, generalizing a result of \cite{CMR-2004}.
\end{abstract}

\section{Introduction}

\subsection{}

A \textit{continuum} is a compact connected metric space. A continuum is called \textit{decomposable} if it is the union of two proper subcontinua; otherwise it is \textit{indecomposable}.

\vspace{2mm}

Throughout, $S^2$ will denote the $2$-dimensional sphere. For a complex polynomial $p : S^2 \to S^2$ of degree $d \geq 2$ all the critical points of which have bounded orbits, the basin of attraction of infinity is a simply connected, fully-invariant region, whose boundary is $J(p)$, the \textit{Julia set} of $p$. Note that $J(p)$ is a fully-invariant continuum. The question of whether the Julia set of a complex polynomial can be an indecomposable continuum has been open for a long time. In \cite{CMR-2004},  Childers, Mayer and Rogers proved a necessary and sufficient condition for $J(p)$ to be indecomposable in terms of the prime ends of its unbounded complementary region in the sphere.

\vspace{2mm}

In this paper we will work with a more general class of branched coverings of the sphere called \textit{polynomial-like} mappings. We say that a branched covering $f : S^2 \to S^2$ of degree $d$ is a \textit{polynomial-like} mapping when there exists a fully invariant simply connected region $R(f)$, whose complement consists of more than one point, containing a critical fixed point with local degree $d-1$. The boundary of $R(f)$ is a fully invariant continuum, which by analogy to the polynomial case we call the \textit{Julia set}\footnote{Note that $R(f)$ is not uniquely defined for a given polynomial-like mapping $f$. Nonetheless, we will refer by $R(f)$ to any region satisfying the conditions in the definition, and, since $R(f)$ will typically be fixed, we will denote its boundary by $J(f)$.} of $f$ and denote it by $J(f)$.

\vspace{2mm}

In \cite{IPRX 22} an example of a family of degree $2$ branched coverings of the sphere having a fully-invariant indecomposable continuum was constructed. However, as far as the author knows, there are currently no known examples of polynomial-like mappings of degree $d \geq 2$ having an indecomposable Julia set. Aiming to shed light on the question of whether such examples exist, our main result is the following.

\begin{theorem}
\label{teo J}
    Let $f : S^2 \to S^2$ be a polynomial-like mapping of degree $d \geq 2$. Then its Julia set $J(f)$ is not the union of finitely many proper indecomposable subcontinua. 
\end{theorem}

\begin{corollary}
\label{coro caraterizacion con prime ends}
    Let $f : S^2 \to S^2$ be a polynomial-like mapping of degree $d \geq 2$. Then $J(f)$ is indecomposable if and only if the impression of some prime end of $R(f)$ is equal to $J(f)$.
\end{corollary}

The polynomial cases of both Theorem \ref{teo J} and Corollary \ref{coro caraterizacion con prime ends} were proved in \cite{MR-93} and \cite{CMR-2004} respectively. However, the strategy implemented in our proof is fundamentally different, so we also provide a new proof of these results for polynomials.

\subsection{Outline of the paper}
In the next section we introduce the main definitions and notation, and recall the minimum requirements to follow the proof of Theorem \ref{teo J}. In Section 3 we prove Theorem \ref{teo J}. The main lemmas used in Section 3 are proved in Section 5. Section 4 is devoted to stating and proving a result needed in Section 5. In the last section we recall the definition of a prime end and its impression as well as show how Corollary \ref{coro caraterizacion con prime ends} follows from Theorem \ref{teo J}.

\subsection{Acknowledgements}
This work originated from the author's undergraduate thesis, done under the advice of Álvaro Rovella, and it would not have existed without him.

Rafael Potrie and Santiago Martinchich generously read preliminary versions of this paper and their help and insightful comments were invaluable.

The author was partially supported by CSIC and a CAP scholarship.

\section{Preliminaries}

\subsection{Indecomposable continua}
A metric space $X$ is called a \textit{continuum} if it is compact and connected. A \textit{subcontinuum} of $X$ is a continuum contained in $X$. We say that $X$ is \textit{decomposable} when it is the union of two proper subcontinua, otherwise we say that $X$ is \textit{indecomposable}.

The main properties of indecomposable continua will be stated later, when needed.

\vspace{1mm}

Given any continuum $X$, an \textit{irreducible decomposition} of $X$ is a finite collection $\{A_1, \ldots, A_n \}$ of indecomposable subcontinua of $X$ such that $X = \bigcup_{i=1}^n A_i$ and none of the $A_j$'s is contained in the union of the others.

\vspace{1mm}

For a continuum $X$ and a subcontinuum $Y \subset X$ we say that a set $A \subset X$ \textit{separates} $Y$ if $Y \setminus A$ is disconnected. If $Z$ is a subset of $Y$, we say that $A$ separates $Y$ \textit{between two points of $Z$} if it separates $Y$ and there are points of $Z$ in distinct connected components of $Y \setminus A$. Since we will be mostly working with subsets of the sphere $S^2$, we will simply say that a set $A \subset S^2$ is \textit{separating} when $A$ separates $S^2$, and will otherwise say that $A$ is \textit{nonseparating}.

\subsection{Polynomial-like mappings}

Here we present the definition and some basic features of the family of maps with which we will be working.

\vspace{2mm}

A map $f: S^2 \to S^2$ is a \textit{branched covering} if it is a covering map up to a finite set of points called \textit{critical} where it is locally conjugate\footnote{i.e. For every critical point $c$, there exist topological disks $D_1$ and $D_2$, with $c \in D_1$, and homeomorphisms $h_1: D_1 \to \D$ and $h_2: \D \to D_2$ such that $f|_{D_1} = h_2 \circ g \circ h_1$, where $g$ denotes the map $\D \to \D / z \mapsto z^k$. } to the map $\D \to \D \ / \ z \mapsto z^k$ for some $k \geq 2$. The number $k - 1$ is the \textit{local degree} of the critical point. The set of critical points of a branched covering $f$ will be denoted by $\textrm{Crit}(f)$. The image $f(c)$ of a critical point $c \in \textrm{Crit}(f)$ is called a \textit{critical value} of $f$. An important property to keep in mind is that removing every critical value of $f$ and its preimage gives us a true covering map (see \cite[Lemma A.11]{BM}).

\begin{remark}
\label{rmk cubrimiento a menos de finitos puntos}
    If $f: S^2 \to S^2$ is a branched covering, there is a finite set $F \subset S^2$ such that the restriction of $f$ to $S^2 \setminus F$ is a covering map onto $S^2 \setminus f(F)$.
\end{remark}

\vspace{2mm}

Given a branched covering $f : S^2 \to S^2$, it follows from the Riemann-Hurwitz formula (see \cite{hurwitz}) that the sum of the local degrees of all critical points of $f$ equals $2d - 2$. Additionally, if a simply connected region $R$ is fully invariant under $f$, then \[ \sum_{c \in \textrm{Crit}(f) \cap R} \textrm{deg}(c) = d - 1, \] where $\textrm{deg}$ denotes the local degree.

\vspace{2mm}

A degree $d$ branched covering $f : S^2 \to S^2$ is called a \textit{polynomial-like mapping} provided that there is a fully-invariant simply connected region $R(f)$ whose complement consists of more than one point, and that there is a critical point $c_\infty \in R(f)$ with local degree $d - 1$.

\vspace{2mm}

The \textit{Julia set} of a polynomial-like mapping $f$ is defined as $J(f) = \partial R(f)$, and it is a fully-invariant continuum. The complement of $R(f)$ in the sphere is the \textit{filled Julia set} of $f$, and we will denote it by $\widehat J(f)$. Note that from the last formula it follows that there are no critical points of $f$ in $R(f)$ other than $c_\infty$.

\section{Proof of Theorem \ref{teo J}}

The proof of Theorem \ref{teo J} will use the next two lemmas, which will be proved in Section 5. Both lemmas will be derived from another key property, namely, that if $\{ A_1, \ldots, A_n \}$ is an irreducible decomposition of $J(f)$, then $f$ acts as a permutation on the set $\{ A_1, \ldots, A_n \}$. This property is obtained from a result regarding \textit{internal composants}, which we will introduce in the next section.

\begin{lemma}
\label{lema preimagen de B-A es conexo}
    Let $f : S^2 \to S^2$ be a polynomial-like mapping of degree $d \geq 2$. If $\{ A_1 , \ldots, A_n \}$ is an irreducible decomposition of $J(f)$, then for every pair $i, j$ the set $f^{-1}(A_i \setminus A_j)$ is connected.
\end{lemma}

\begin{lemma}
\label{lema L}
    Let $f : S^2 \to S^2$ be a polynomial-like mapping of degree $d \geq 2$, and let $\alpha$ be a simple curve joining $f(c_\infty)$ with a point $y \in J(f)$, such that $\alpha \setminus \{ y \}$ is contained in $R(f)$. Suppose that $\{ A_1, \ldots, A_n \}$ is an irreducible decomposition of $J(f)$ and let $L = \alpha \cup A_j$, where $y \in A_j$. Then $f^{-1}(L)$ separates $S^2$ and every point in $J(f) \setminus A_j$ has preimages in at least two distinct connected components of $S^2 \setminus f^{-1}(L)$. 
\end{lemma}

\vspace{2mm}

We can now prove Theorem \ref{teo J}.

\vspace{1mm}

\textit{Proof of Theorem \ref{teo J}.}
    Recall that we want to show that $J(f)$ is not the union of a finite number of indecomposable proper subcontinua. Equivalently, that there are no irreducible decompositions of $J(f)$ with more than one element.

    \vspace{2mm}

    Let us suppose by contradiction that $\{ A_1, \ldots, A_n \}$ is an irreducible decomposition of $J(f)$ with $n \geq 2$.

    \vspace{2mm}

    Let $\alpha$ be a curve as described in Lemma \ref{lema L}, joining $f(c_\infty)$ with some point $y \in A_j$, and let $L = \alpha \cup A_j$. If $i \neq j$, then $A_i \setminus A_j$ is contained in $J(f) \setminus A_j$, so Lemma \ref{lema L} ensures that every point in $A_i\setminus A_j$ has preimages in two connected components of the complement of $f^{-1}(L)$. Since $L$ is disjoint from $A_i \setminus A_j$, then $f^{-1}(A_i \setminus A_j)$ must be disjoint from $f^{-1}(L)$. We conclude that $f^{-1}(A_i \setminus A_j)$ is  not connected, contradicting Lemma \ref{lema preimagen de B-A es conexo}
    
\qed

\section{A theorem regarding internal composants}

The purpose of this section is to prove Theorem \ref{teo K}, which will be used in the next section, where we prove Lemmas \ref{lema preimagen de B-A es conexo} and \ref{lema L}. Theorem \ref{teo K} deals with the notion of \textit{internal composants} of an indecomposable continuum, first introduced by Krasinkiewicz in \cite{krasi 72}. 

It should be mentioned that the only parts of this section that will be used later are contained in Subsection 4.1, so the reader can skip the rest of the section and have no problem following the arguments in future sections.

\subsection{Composants}

The \textit{composant} of a point $p$ in a continuum $X$ is defined as the union of all proper subcontinua of $X$ containing $p$. The following facts can be found in \cite[Ch. 3-8]{HY}

\begin{proposition}
\label{prop compos es union numerable de continuos y es densa}
    Every composant of a continuum $X$ is dense in $X$ and it is the union of a countable family of subcontinua of $X$.
\end{proposition}

\begin{proposition}
\label{prop en indesc. subcontinuos no tienen interior}
    A continuum is indecomposable if and only if all of its subcontinua have empty interior.
\end{proposition}

It can happen that two distinct composants of a continuum overlap, as is the case with some decomposable continua such as a line segment. However, for indecomposable continua the following holds.

\begin{proposition}
\label{prop compos en indescomponible son disjuntas}
    Two composants of an indecomposable continuum are either equal or disjoint.
\end{proposition}

A composant $C$ of a continuum $X$ is said to be \textit{internal} if every continuum $L$ intersecting $C$ but not contained in $X$ intersects every composant of $X$. Otherwise $C$ is said to be \textit{external}. 

With these concepts in mind we can state the main result of this section. This statement coincides with \cite[Thm 3.1]{MR-93}. Here we give a proof based on [krasi74] for completeness.

\begin{theorem}
\label{teo K}
    Let $X \subset S^2$ be a nonseparating continuum, and let $Y$ be an indecomposable continuum contained in the boundary of $X$. If there exists a subcontinuum $Z$ of $X$ intersecting an internal composant of $Y$, then either $Z$ is contained in $Y$ or $Y$ is contained in $Z$.
\end{theorem}

\subsection{Proof of Theorem \ref{teo K}}

\begin{figure}[ht]
    \centering
    \def\svgwidth{0.4\columnwidth}
    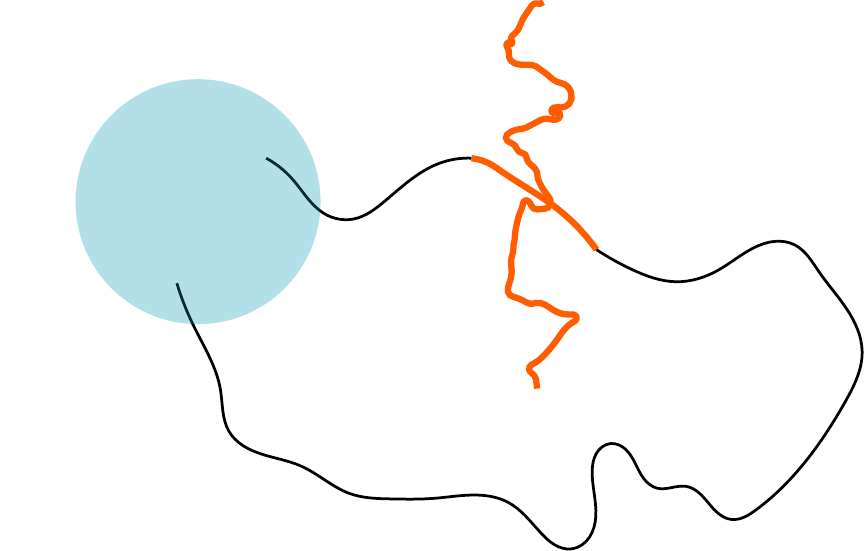

    \caption{An illustration of Theorem \ref{teo krasi}.}
    \label{krasi}        
\end{figure}

The proof will be based on the following result due to Krasinkiewicz (see Figure \ref{krasi} for a picture).

\begin{theorem}
\label{teo krasi}
    Let $Y \subset S^2$ be an indecomposable continuum and let $L$ be a continuum intersecting an internal composant $C$ of $Y$. If $L$ is not contained in $Y$ and $Y$ is not contained in $L$, then for every disk $D$ disjoint from $L$ and intersecting $Y$ there is a continuum $A \subset C$ such that $A \cup D$ separates $S^2$ between two points of $L$.
\end{theorem}

\proof
    In \cite[2.3]{krasi 74}, take $X = Y$ and $K = D$, and let $R$ be the collection of all composants of $Y$.
    
\qed

\vspace{2mm}

The first step will be to prove a weakened version of the theorem, where we suppose that $Z$ is nonseparating.

\begin{proposition}
\label{prop teo K weak}
    Let $X \subset S^2$ be a nonseparating continuum and let $Y$ be an indecomposable continuum contained in the boundary of $X$. If $Z$ is a subcontinuum of $X$ intersecting an internal composant of $Y$ and $Z$ does not separate $S^2$, then either $Y \subset Z$ or $Z \subset Y$.
\end{proposition}

To prove this proposition we will need some lemmas.

\begin{lemma}
\label{lema si S seprara continuo, las cc tocan S}
    If $A$ is a continuum and $S \subset A$ is a nonempty compact subset, then, for every connected component $B$ of $A \setminus S$, the boundary of $B$ intersects $S$.
\end{lemma}

\proof
    Recall that the boundary bumping lemma (\cite[\textsection 47, III, Thm. 3]{kuratowski}) states that, if $Y$ is a closed subset of a continuum $X$, then every connected component of $Y$ meets the boundary of $Y$.

    \vspace{2mm}

    Let $x$ be any point in $A \setminus S$ and let $B$ be the connected component of $x$ in $A \setminus S$. For every $n \geq 1$, let $U_n$ be the open $\tfrac{1}{n}$-neighborhood of $S$. Suppose $n$ is large enough so that $x \notin U_n$. Since $A \setminus U_n$ is closed in $A$, by boundary bumping, the connected component $C_n$ of $x$ in $A \setminus U_n$ intersects the boundary of $U_n$. But $C_n$ is contained in $B$, meaning that, for $n$ arbitrarily large, there is a point $q_n \in B$ at distance $1/n$ from $S$. By compactness of $S$, this implies that $\overline{B}$ intersects $S$.
    
\qed

\vspace{2mm}

\begin{lemma}
\label{lema C inter Z es conexo}
Let $X \subset S^2$ be a nonseparating continuum, and let $Y$ be an indecomposable continuum contained in the boundary of $X$. If $Z$ is a nonseparating subcontinuum of $X$ intersecting an internal composant $C$ of $Y$, then $C \cap Z$ is connected.  
\end{lemma}

\proof
    Let us suppose by contradiction that $C \cap Z$ is not connected, and take points $x$ and $y$ in distinct connected components of $C \cap Z$. Let $K \subset C$ be a continuum containing $x$ and $y$. Then $K \cap Z$ must be disconnected, so, by Theorem \ref{teo janis} we have that $K \cup Z$ separates $S^2$. Because $X$ is nonseparating and $K \cup Z$ is contained in $X$, some connected component $U$ of $S^2 \setminus (K \cup Z)$ must be contained in $X$. We have that the boundary of $U$ is contained in $K \cup Z$, but it is not fully contained in $Z$, since $Z$ is nonseparating.

    \vspace{2mm}

    Take a point $c \in \partial U \setminus Z$ and a neighborhood $V$ of $c$ disjoint from $Z$. For some $u \in U \cap V$ let $\alpha$ be a curve contained in $V \cap \overline{U}$ joining $u$ with a point $p \in \partial U$, such that $\alpha \setminus \{ p \} \subset U$. Then $p \in K$. But, because $U$ is an open set contained in $X$ and $Y \subset \partial X$, we have that $U$ does not intersect $Y$. This means that no point of $\alpha$ other than $p$ lies in $Y$, contradicting the fact that $C$ is internal.
    
\qed

\vspace{2mm}

\begin{lemma}
\label{lema C inter Z es compacto}
    Let $A$ and $B$ be subcontinua of $S^2$ such that, for some composant $C$ of $A$, the intersection $C \cap B$ is connected. Then either $A \subset B$ or $C \cap B$ is compact. 
\end{lemma}

\proof
    Suppose that $\overline{C \cap B} = A$. Then $A \subset \overline C \cap B = A \cap B$, meaning that $A \subset B$.

    \vspace{2mm}
    
    Alternatively, $\overline{C \cap B}$ is strictly contained in $A$. Since $\overline{C \cap B}$ is connected, it is a proper subcontinuum of $A$ intersecting $C$. By Proposition \ref{prop compos en indescomponible son disjuntas}, this implies that $\overline{C \cap B}$ is contained in $C$. It follows that $\overline{C \cap B} \subset C \cap B$, which means that $C \cap B$ is compact.
    
\qed

\proof[Proof of Proposition \ref{prop teo K weak}.]
    Let $C$ be an internal composant of $Y$ intersecting $Z$. From Lemmas \ref{lema C inter Z es conexo} and \ref{lema C inter Z es compacto} we get that $Z \cap C$ is connected and either $Y$ is contained in $Z$ or $C \cap Z$ is compact.

    Suppose by contradiction that $Y$ is not contained in $Z$ and $Z$ is not contained in $Y$. Take any point $z_0 \in Z \setminus Y$ and let $Z_0$ be the connected component of $Z \setminus C$ containing $z_0$. Because $C \cap Z$ is compact, we know, via Lemma \ref{lema si S seprara continuo, las cc tocan S}, that the closure $\overline{Z_0}$ intersects $C$. So $\overline{Z_0}$ is a continuum intersecting an internal composant of $Y$ and it is not contained in nor contains $Y$. Hence, given any disk $D$ intersecting $Y$ but not $\overline{Z_0}$, Theorem \ref{teo krasi} guarantees that there exists a continuum $A \subset C$ such that $A \cup D$ separates $\overline{Z_0}$. But $D \cap \overline{Z_0} = \emptyset$, so $A$ must separate $\overline{Z_0}$. Since $A$ is compact, this means that $A$ intersects $Z_0$, so $C \cap Z_0 \neq \emptyset$, contradicting the definition of $Z_0$.
    
\qed

\vspace{2mm}

From this we can prove the general version of Theorem \ref{teo K}.

\vspace{2mm}

\textit{Proof of Theorem \ref{teo K}.}
    Let $\widehat Z$ be the union of $Z$ with all the connected components of its complement in the sphere which are contained in $X$. That is, \[ \widehat Z  = Z \cup \bigcup \{ U: U \textrm{ is a connected component of  $S^2 \setminus Z$ with $U \subset X$} \}.\] Since $X$ is nonseparating, so is $\widehat Z$. From Proposition \ref{prop teo K weak} it follows that $Y \subset \widehat Z$. But $\widehat Z \setminus Z$ is contained in the interior of $X$, while $Y$ is contained within the boundary of $X$. Hence, $Z$ contains $Y$.
    
\qed

\section{Indecomposable subcontinua of Julia sets}

The goal of this section is to prove Lemmas \ref{lema preimagen de B-A es conexo} and \ref{lema L}.

\subsection{Invariance of the decomposition}

The following is an important result, which will be key to prove Lemmas \ref{lema preimagen de B-A es conexo} and \ref{lema L}.

\begin{lemma}
\label{lema f permuta los indescomponibles}
    Let $f : S^2 \to S^2$ be a polynomial-like mapping of degree $d \geq 2$. If $\{ A_1 , \ldots, A_n \}$ is an irreducible decomposition of $J(f)$, then for every $j$ there exists exactly one $k$ such that $f^{-1}(A_j) = A_k$.
\end{lemma}

In other words, Lemma \ref{lema f permuta los indescomponibles} states that $f$ acts as a permutation on the set $\{ A_1 , \ldots, A_n\}$, meaning that there is some iterate of $f$ for which the $A_i$'s are fully invariant.

\begin{proposition}
\label{prop compos de indescomponibles no tienen interior}
    Every composant of an indecomposable continuum $X$ has empty interior in $X$.
\end{proposition}

\proof

This follows from Propositions \ref{prop en indesc. subcontinuos no tienen interior} and \ref{prop compos es union numerable de continuos y es densa} and Baire's category theorem.

\qed

\begin{lemma}
\label{lema imagen de compos no tiene interior en J}
    Let $f : S^2 \to S^2$ be a poynomial-like mapping and suppose that $A$ is an indecomposable subcontinuum of $J(f)$. Then, the image $f(C)$ of any composant $C$ of $A$ has empty interior in $J(f)$.
\end{lemma}

\proof
    Suppose by contradiction that $U\cap J(f)$ is contained in $f(C)$ for some open set $U$ intersecting $J(f)$. Since $J(f)$ has no isolated points and $f$ has only finitely many critical values, without loss of generality we can assume that $U$ contains no critical values. Moreover, we can assume that $U$ is well covered, that is, such that the preimage of $U$ is a disjoint union of open sets $U_1, \ldots, U_d$, for which each restriction $f|_{U_i} : U_i \to U$ is a homeomorphism. Hence, for every $i$ we have that $f$ maps subsets of $U_i$ with empty interior in $J(f)$ to subsets of $U$ with empty interior in $J(f)$. Note that $U \cap J(f)$ is contained in $\bigcup_{i = 1}^d f(U_i \cap C)$.

    \vspace{1mm}

    By Proposition \ref{prop compos es union numerable de continuos y es densa} we have that $C$ is a countable union of proper subcontinua of $A$. Then, for each $i$ the intersection $U_i \cap C$ is a countable union of compact sets. We know that $C$ has empty interior in $A$ thanks to Proposition \ref{prop compos de indescomponibles no tienen interior}. Therefore, every $U_i \cap C$ has empty interior in $J(f)$.
    
    By the above we have that $f(U_i \cap C)$ is a countable union of compact sets with empty interior in $J(f)$. It follows from Baire's category theorem that $\bigcup_{i = 1}^d f(U_i \cap C)$ has empty interior in $J(f)$, contradicting the fact that $U \cap J(f)$ is contained in $\bigcup_{i = 1}^d f(U_i \cap C)$.

\qed

\vspace{3mm}

Suppose that $\{A_1, \ldots, A_n \}$ is an irreducible decomposition of a continuum $X$. If some $A_k$ had no interior in $X$, then the union of all the others would be dense in $X$ and, therefore, equal to $X$, meaning that $A_j$ would be contained in the union of the others. We obtain the following.

\begin{remark}
\label{rmk Ai tienen interior en J}
    If $f : S^2 \to S^2$ is a polynomial-like mapping of degree $d \geq 2$ and $\{ A_1, \ldots, A_n \}$ is an irreducible decomposition of $J(f)$, then every $A_i$ has nonempty interior in $J(f)$.
\end{remark}

\begin{lemma}
\label{lema compos internas no se cortan}
    Let $f : S^2 \to S^2$ be a polynomial-like mapping of degree $d \geq 2$. If $\{ A_1 , \ldots, A_n \}$ is an irreducible decomposition of $J(f)$, then for every $j \neq k$, no internal composant of $A_k$ intersects $A_j$.
\end{lemma}

\proof
Let us suppose by contradiction that $A_j$ intersects some internal composant $C_k$ of $A_k$, with $j \neq k$, and let $X = \widehat J(f)$, $Y = A_k$ and $Z = A_j$. We have that both $Y$ and $Z$ are subcontinua of $J(f)$, which is the boundary of $X$. Additionally, $A_j$ is not contained in $A_k$, so, thanks to Theorem \ref{teo K}, we get that $A_k$ is contained in $A_j$. Thus, the decomposition is not irreducible.

\qed

\vspace{1mm}

\begin{lemma}
\label{lema f(Aj) contenido en Ak}
    Let $f : S^2 \to S^2$ be a polynomial-like mapping of degree $d \geq 2$. If $\{ A_1 , \ldots, A_n \}$ is an irreducible decomposition of $J(f)$, then for each pair $j, k$, if a point $x \in A_j$ is such that $f(x)$ lies in an internal composant of $A_k$, then  $f(A_j)$ is contained in $A_k$.
\end{lemma}

\proof
    Let $C^j_x$ be the composant of $x$ in $A_j$ and let $C^k_{f(x)}$ be the (internal) composant of $f(x)$ in $A_k$. We will see that $f(C_x^j)$ is contained in $A_k$. Since $C_x^j$ is dense in $A_j$, this will finish the proof.

    \vspace{2mm}
    
    Given any proper subcontinuum $K$ of $A_j$ containing $x$, define $Z = f(K)$, $Y = A_k$ and $X = \widehat J(f)$. We have that $X$, $Y$ and $Z$ satisfy the hypothesis of Theorem \ref{teo K}, therefore either $A_k$ is contained in $f(K)$ or $f(K)$ is contained in $A_k$. By Lemma \ref{lema imagen de compos no tiene interior en J} we know that $f(K)$ has empty interior in $J(f)$, while $A_k$ does not, as stated in Remark \ref{rmk Ai tienen interior en J}, making it impossible for $A_k$ to be contained in $f(K)$. Hence, $f(K)$ must be contained in $A_k$.

    \vspace{2mm}
    
    We have proved that the image of every proper subcontinuum of $A_j$ containing $x$ is contained in $A_k$, thus $f(C^j_x)$ is also contained in $A_k$, as we wanted.
    
\qed

\vspace{3mm}

In \cite{krasi 72} Krasinkiewicz proved that the union of all external composants of an indecomposable continuum $X$ is a first category subset of $X$. In particular, we will use the fact that every indecomposable continuum has some internal composant, so we highlight it as a remark.

\begin{remark}
\label{rmk existen compos internas}
    Every indecomposable continuum has an internal composant.
\end{remark}

\vspace{1mm}

\begin{lemma}
\label{lema preimagen de Aj es union de Aks}
    Let $f : S^2 \to S^2$ be a polynomial-like mapping of degree $d \geq 2$ and suppose that $\{ A_1 , \ldots, A_n \}$ is an irreducible decomposition of $J(f)$. For every $j$ define $H_j$ as \[ H_j = \{ k: f(A_k) \textrm{ intersects an internal composant of $A_j$} \}. \] Then $f^{-1}(A_j) = \bigcup_{k \in H_j} A_k$.
\end{lemma}

\proof
    The inclusion $\bigcup_{k \in H_j} A_k \subset f^{-1}(A_j)$ follows from Lemma \ref{lema f(Aj) contenido en Ak}. To prove the converse, take $x \in f^{-1}(A_j)$ and let us suppose by contradiction that $x$ does not belong to any $A_k$ with $k \in H_j$. By Remark \ref{rmk existen compos internas} there exists some internal composant $C_j$ of $A_j$. Let $V$ be an open neighborhood of $x$ disjoint from $\bigcup_{k \in H_j} A_k$. Because $f$ is open, we know that $f(V)$ is an open neighborhood of $f(x) \in A_j$, so $C_j$ intersects $f(V)$. In other words, there is some $z \in V$ such that $f(z)$ lies in an internal composant of $A_j$, contradicting the definition of $H_j$.
    
\qed

\vspace{2mm}

As a corollary we obtain Lemma \ref{lema f permuta los indescomponibles}.

\proof[Proof of Lemma \ref{lema f permuta los indescomponibles}]

\vspace{2mm}

    Note that from our last lemma it follows that $H_j$ is nonempty for every $j$. To finish the proof it suffices to show that $H_j$ is disjoint from $H_k$ if $j \neq k$. If there were some $ \ell \in H_j \cap H_k$, then, by Lemma \ref{lema f(Aj) contenido en Ak}, we would have that $f(A_\ell) \subset A_j \cap A_k$. However, we would also have that $f(A_\ell)$ intersects some internal composant of $A_j$, meaning that $A_k$ intersects an internal composant of $A_j$, contradicting Lemma \ref{lema compos internas no se cortan}.

\qed

\subsection{Proof of Lemma \ref{lema preimagen de B-A es conexo}}

\proof[Proof of Lemma \ref{lema preimagen de B-A es conexo}.]

    We can assume $i \neq j$. By Lemma \ref{lema f permuta los indescomponibles}, for some $k \neq \ell$ we have \[ f^{-1}(A_i \setminus A_j) = f^{-1}(A_i) \setminus f^{-1}(A_j) = A_k \setminus A_\ell \]
    
    Remark \ref{rmk existen compos internas} assures that there exists an internal composant $C$ of $A_k$. From Lemma \ref{lema compos internas no se cortan} we know that $C \subset A_k \setminus A_\ell$. Since $C$ is connected and dense in $A_k$, we have that $A_k \setminus A_\ell$ is connected.

\qed

\subsection{Proof of Lemma \ref{lema L}}

First we will prove the next partial result.

\begin{lemma}
\label{lema preimagen separa}
    Under the hypotheses of Lemma \ref{lema L}, suppose that $K$ is a subcontinuum of the filled Julia set $\widehat J(f)$, with $f^{-1}(y) \subset K$. Then $f^{-1}(\alpha) \cup K$ separates $S^2$.
\end{lemma}

The proof of Lemma \ref{lema preimagen separa} will employ the following theorem by Janiszewski (see \cite{janis}).

\begin{theorem}
\label{teo janis}
    Let $A$ and $B$ be two proper subcontinua of the sphere. Then
    \begin{enumerate}
        \item If $A \cap B$ is connected and $A$ and $B$ are nonseparating, then $A \cup B$ is nonseparating.
        \item If $A \cap B$ is disconnected, then $A \cup B$ separates $S^2$.
    \end{enumerate}
\end{theorem}

\vspace{1mm}

\proof[Proof of Lemma \ref{lema preimagen separa}]
    Because $c_\infty$ has local degree $d-1$, there are $d$ lifts of $\alpha$, each going from $c_\infty$ to some point in the preimage of $y$. Moreover, since $\alpha$ is simple and contains no critical values other than $f(c_\infty)$ (and possibly $y$), the intersection of two distinct lifts of $\alpha$ is contained\footnote{In fact, if $y$ is chosen to be a regular value of $f$, then the intersection of two distinct lifts of $\alpha$ will be exactly $\{c_\infty\}$.} in $\{c_\infty\} \cup f^{-1}(y)$.

    \begin{figure}[ht]
        \centering
    \def\svgwidth{.92\columnwidth}
    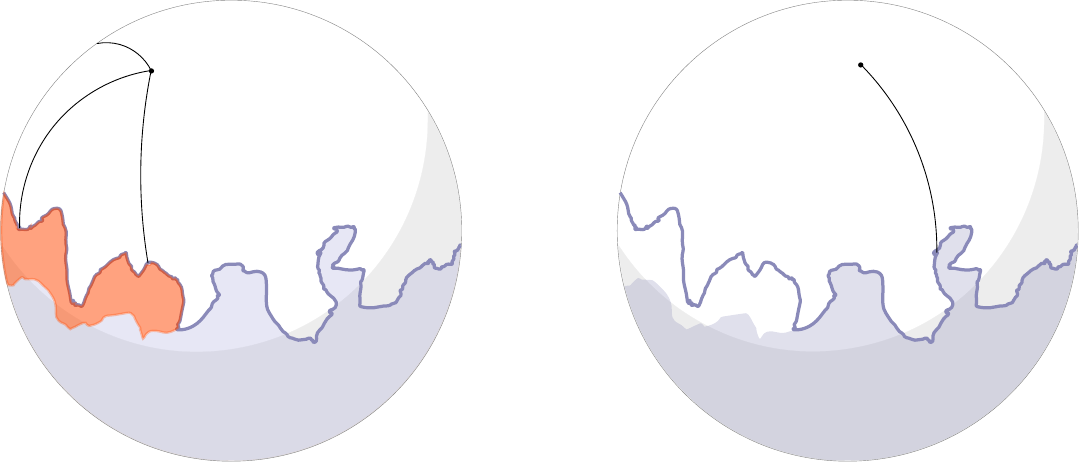

        \caption{}
        \label{lifts}
    \end{figure}
    
    Let $\hat \alpha_1$ and $\hat \alpha_2$ be two distinct lifts of $\alpha$, and let $\hat y_1, \hat y_2 \in f^{-1}(y)$ be their respective endpoints (as in Figure \ref{lifts}). Because $\hat y_1, \hat y_2 \in K$, the sets $K_1 = K \cup \hat \alpha_1$ and $K_2 = K \cup \hat \alpha_2$ are connected, and their intersection is $K_1 \cap K_2 = \{ c_\infty \} \cup K$, which is disconnected since $K$ is contained in $\widehat J(f)$. By Theorem \ref{teo janis}, this implies that $K_1 \cup K_2$ separates $S^2$. But, since $K_1 \cup K_2$ is compact and $(K \cup f^{-1}(\alpha)) \setminus (K_1 \cup K_2)$ is contained in $f^{-1}(\alpha)$, which has empty interior, we have that $K \cup f^{-1}(\alpha)$ separates $S^2$.
    
\qed

\begin{lemma}
\label{lema imagen de cc es cc}
    Let $f: S^2 \to S^2$ be a continuous, open map and $L \subset S^2$ a compact set. If $U$ is a connected component of $S^2 \setminus f^{-1}(L)$, then $f(U)$ is a connected component of $S^2 \setminus L$.
\end{lemma}

\proof
    Since $f(U)$ is a connected set contained in $S^2 \setminus L$, it is enough to show that it is both open and closed as a subset of $S^2 \setminus L$. 
    
    Note that $S^2 \setminus L$ is open in $S^2$. We have that $f(U)$ is open in $S^2$, because $f$ is open, so $f(U)$ is open in $S^2 \setminus L$. 

    On the other hand, we have $\partial f(U) \subset f(\partial U)$. Indeed, given a point $x \in \partial f(U)$, there is a sequence $\{u_n\}_n$ of points in $U$ such that $f(u_n) \xrightarrow[n]{} x$. Let $u \in \overline{U}$ be a limit point of $\{u_n\}_n$. Then $f(u) = x$. So, because $f(U)$ is open and $f(u) \in \partial f(U)$, it cannot be $u \in U$. Thus, $u \in \partial U$, so $x \in f(\partial U)$.

    By the above we get that $\partial f(U) \subset f(\partial U) \subset L$, which means that $f(U) = \overline{f(U)} \cap (S^2 \setminus L)$. So $f(U)$ is a closed subset of $S^2 \setminus L$.
\qed

\vspace{2mm}

\begin{corollary}
\label{coro hay preimagen en al menos dos cc}
    If $f: S^2 \to S^2$ is a continuous, open map and $L \subset S^2$ is a nonseparating compact set, then every point in $S^2 \setminus L$ has at least one preimage in every connected component of $S^2 \setminus f^{-1}(L)$.
\end{corollary}

\proof
    By the above lemma, given any connected component $U$ of $S^2 \setminus f^{-1}(L)$, we have that $f$ maps $U$ to a whole connected component of $S^2 \setminus L$. Because $L$ is nonseparating, this means that $f(U) = S^2 \setminus L$. 

\qed

\begin{lemma}
\label{lema rellenos}
    Let $f: S^2 \to S^2$ be a polynomial-like mapping of degree $d \geq 2$ and suppose that $\{ A_1, \ldots, A_n \}$ is an irreducible decomposition of $J(f)$. For every $i =1, \ldots, n$ denote by $U_i$ the connected component of $S^2 \setminus A_i$ containing $R(f)$ and consider the `filled' set $\widehat A_i = S^2 \setminus U_i$. Then
    \begin{enumerate}
        \item $\widehat A_i$ is a nonseparating subcontinuum of $\widehat J(f)$ and
        \item if $f^{-1}(A_j) = A_k$, then $f^{-1}(\widehat A_j) = \widehat A_k$.
    \end{enumerate}
\end{lemma}

\proof
    First note that \[ \widehat A_i =  A_i \cup \bigcup \{ U : U \textrm{ connected component of $ S^2 \setminus A_i$ with $U \neq U_i$} \} \]
    \begin{figure}[ht]
        \centering
    \def\svgwidth{0.4\columnwidth}
    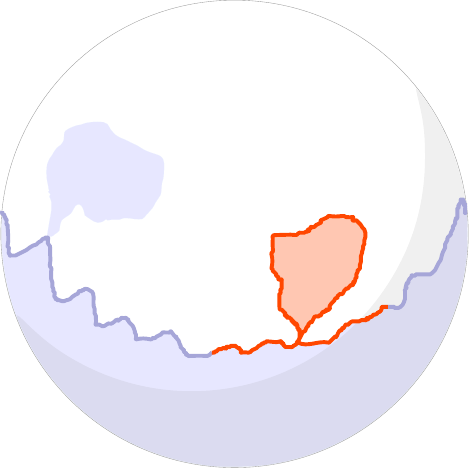

        \caption{$A_i$ can separate $S^2$ in many connected components.}
        \label{rellenos}
    \end{figure}
    and we have that the boundary of any connected component $U$ of $S^2 \setminus A_i$ is contained in $A_i$, since $A_i$ is compact. This means that \[ \widehat A_i = A_i \cup \bigcup \{ \overline U : U \textrm{ connected component of $ S^2 \setminus A_i$ with $U \neq U_i$} \}. \] So $\widehat A_i$ is connected, as it is the union of connected sets, all of which intersect $A_i$.

    Additionally, because $U_i$ is open, connected and contains $R(f)$, we have that $\widehat A_i$ is compact, nonseparating and contained in $\widehat J(f)$. This proves claim (1).

    \vspace{1mm}

    Now suppose that $f^{-1}(A_j) = A_k$. By Lemma \ref{lema imagen de cc es cc}, for any connected component $U$ of $S^2 \setminus A_k$ it holds that $f(U)$ is equal to a connected component of $S^2 \setminus A_j$. Therefore, if $f(U)$ intersects $U_j$, we have that $f(U) = U_j$. Since $R(f)$ is contained in $U$ and it is fully invariant, this means that there are points of $R(f)$ in U, so it must be $U = U_k$. This proves that $f^{-1}(U_j) \subset U_k$. For the converse note that $f(U_k)$ is connected and contains $R(f)$, so we have that $f(U_k)$ is contained in $U_j$. We have shown that $f^{-1}(U_j) = U_k$, and this finishes our proof, since \[ f^{-1}(\widehat A_j) = S^2 \setminus f^{-1}(U_j) = S^2 \setminus U_k = \widehat A_k. \] 
\qed

\vspace{2mm}

\proof[Proof of Lemma \ref{lema L}]
    We have, due to Lemma \ref{lema f permuta los indescomponibles}, that $f^{-1}(A_j) = A_k$ for some $k$. 
    
    For each $i = 1, \ldots, n$, let $\widehat A_i$ be as in Lemma \ref{lema rellenos}. By Lemma \ref{lema rellenos} we have that $\widehat A_i$ is a nonseparating subcontinuum of $\widehat J(f)$ for every $i$, and $f^{-1}(\widehat A_j) = \widehat A_k$. Since $y \in A_j$, we know that $f^{-1}(y) \subset A_k \subset \widehat{A}_k$. Therefore, by Lemma \ref{lema preimagen separa} both $f^{-1}(\alpha) \cup A_k$ and $f^{-1}(\alpha) \cup \widehat{A}_k$ separate the sphere. This proves that $f^{-1}(L)$ separates $S^2$, since  \[f^{-1}(L) = f^{-1}(\alpha \cup A_j) =  f^{-1}(\alpha) \cup A_k. \]

    It remains to show that the preimage of any point in $J(f) \setminus A_j$ intersects at least two distinct connected components of the complement of $f^{-1}(L)$.

    We have established that $f^{-1}(\alpha \cup \widehat A_j) = f^{-1}(\alpha) \cup \widehat A_k$ separates $S^2$. Note that each connected component of $S^2 \setminus f^{-1}(L)$ is either contained in $\widehat A_k$ or disjoint from it, as it does not contain points of the boundary of $\widehat A_k$.  Hence, any connected component of $S^2 \setminus f^{-1}(\alpha \cup \widehat A_j)$ is also a connected component of $S^2 \setminus f^{-1}(L)$.

    Let $x$ be a point in $J(f) \setminus A_j$. Then $x$ cannot belong in $\widehat A_j$, since $J(f)$ is the boundary of $R(f)$. Thus, $x \in S^2 \setminus \widehat A_j$. Because $\alpha \setminus \{ y \}$ is contained in $R(f)$, the intersection $\alpha \cap \widehat A_j$ consists of a single point. Then, since $\alpha$ is simple and $\widehat A_j$ is nonseparating, by Theorem \ref{teo janis} the set $\alpha \cup \widehat A_j$ is nonseparating. Therefore, Corollary \ref{coro hay preimagen en al menos dos cc} applies here, meaning that $f^{-1}(x)$ intersects all connected components of $S^2 \setminus f^{-1}(\alpha \cup \widehat A_j)$. If $U$ and $V$ are two distinct connected components of $S^2 \setminus f^{-1}(\alpha)$ containing points of $f^{-1}(x)$, by the above we have that $U$ and $V$ are also connected components of $S^2 \setminus f^{-1}(L)$, as we wanted. 
    
\qed

\section{Prime ends and indecomposability of the Julia set}

\subsection{Prime ends} We will briefly recall the definition of a prime end and its impression. We refer the reader to \cite{milnor} for further details on the theory of prime ends.

\vspace{1mm}

Let $U \subset S^2$ be a simply connected open set whose complement has more than one point. A \textit{crosscut} of $U$ is a simple arc $A 
\subset \overline U$ intersecting the boundary of $U$ only at its endpoints. A \textit{fundamental chain} in $U$ is a sequence $\{ R_i \}_{i \in \mathbb{N}}$ of nested regions in $U$, such that each $R_i$ is a connected component of $U \setminus A_i$, where $\{A_i\}_{i \in \mathbb{N}}$ is a sequence of crosscuts of $U$ whose diameter tends to $0$. We define two fundamental chains of $U$ to be equivalent when every region of one contains a region of the other and vice versa.

\vspace{1mm}

With the above definitions, a \textit{prime end} of $U$ is an equivalence class of fundamental chains in $U$. if $\varepsilon$ is a prime end of $U$ defined by a fundamental chain $\{ R_i \}_{i \in \mathbb N}$, the \textit{impression} of $\varepsilon$ is the set \[ \textrm{Im}(\varepsilon) = \bigcap_{i} \overline{R_i}. \] It can be shown that the impression of a prime end of $U$ is always a continuum contained in the boundary of $U$ (\cite[Lemma 17.7]{milnor}).

\subsection{Rutt's theorems and Corollary \ref{coro caraterizacion con prime ends}}

In \cite{rutt}, E. Rutt proved the following two results.

\begin{theorem}[Rutt]
\label{teo Rutt 1}
    Let $U \subset S^2$ be a simply connected open region whose boundary $\partial U$ has more than one point. If $\partial U$ is indecomposable, then there exists a prime end of $U$ whose impression is equal to $\partial U$.
\end{theorem}

\begin{theorem}[Rutt]
\label{teo Rutt 2}
    Let $U \subset S^2$ be a simply connected open region whose boundary $\partial U$ has more than one point. If the impression of some prime end of $U$ is the entire $\partial U$, then either $\partial U$ is indecomposable or it is the union of two proper indecomposable subcontinua.
\end{theorem}

\vspace{2mm}

In the same article Rutt also provided examples showing that neither converse is true. However, thanks to Theorem \ref{teo J} we can prove that, in our context, the existence of a prime end with impression equal to the Julia set is an equivalent condition to the Julia set being indecomposable.

\proof[Proof of Corollary \ref{coro caraterizacion con prime ends}]
    Suppose that $J(f)$ is an indecomposable continuum, and let $U$ be a connected component of $J(f)$ with $\partial U = J(f)$, such as $R(f)$. Then, by Theorem \ref{teo Rutt 1}, there is a prime end of $U$ whose impression is $J(f)$.
    Conversely, if the impression of some prime end of a connected component $U$ of $S^2 \setminus J(f)$ is equal to $J(f)$, then the boundary of $U$ is equal to $J(f)$ and, by Theorem \ref{teo Rutt 2}, we have that $J(f)$ is indecomposable or the union of two indecomposable subcontinua. Due to Theorem \ref{teo J} we know that the latter cannot happen, so $J(f)$ must be indecomposable.
\qed

\vspace{2mm}

Note that we have actually proved the following.

\begin{corollary}
    Let $f : S^2 \to S^2$ be a polynomial-like mapping of degree $d \geq 2$. Then $J(f)$ is indecomposable if and only if, for every connected component $U$ of $S^2 \setminus J(f)$ such that $\partial U = J(f)$, there is a prime end of $U$ whose impression is equal to $J(f)$.
\end{corollary}

\end{document}

%% file: Pictures/krasi.pdf_tex
\begingroup%
  \makeatletter%
  \providecommand\color[2][]{%
    \errmessage{(Inkscape) Color is used for the text in Inkscape, but the package 'color.sty' is not loaded}%
    \renewcommand\color[2][]{}%
  }%
  \providecommand\transparent[1]{%
    \errmessage{(Inkscape) Transparency is used (non-zero) for the text in Inkscape, but the package 'transparent.sty' is not loaded}%
    \renewcommand\transparent[1]{}%
  }%
  \providecommand\rotatebox[2]{#2}%
  \newcommand*\fsize{\dimexpr\f@size pt\relax}%
  \newcommand*\lineheight[1]{\fontsize{\fsize}{#1\fsize}\selectfont}%
  \ifx\svgwidth\undefined%
    \setlength{\unitlength}{414.49757766bp}%
    \ifx\svgscale\undefined%
      \relax%
    \else%
      \setlength{\unitlength}{\unitlength * \real{\svgscale}}%
    \fi%
  \else%
    \setlength{\unitlength}{\svgwidth}%
  \fi%
  \global\let\svgwidth\undefined%
  \global\let\svgscale\undefined%
  \makeatother%
  \begin{picture}(1,0.63742099)%
    \lineheight{1}%
    \setlength\tabcolsep{0pt}%
    \put(0,0){\includegraphics[width=\unitlength,page=1]{krasi.pdf}}%
    \put(-0.00320418,0.39943218){\color[rgb]{0,0.42352941,0.50980392}\transparent{0.5}\makebox(0,0)[lt]{\lineheight{10.25}\smash{\begin{tabular}[t]{l}$D$\end{tabular}}}}%
    \put(0.70415042,0.46157701){\color[rgb]{1,0.36470588,0}\makebox(0,0)[lt]{\lineheight{10.25}\smash{\begin{tabular}[t]{l}$L$\end{tabular}}}}%
    \put(0.27096424,0.03021883){\color[rgb]{0,0,0}\makebox(0,0)[lt]{\lineheight{10.25}\smash{\begin{tabular}[t]{l}$C$\end{tabular}}}}%
    \put(0,0){\includegraphics[width=\unitlength,page=2]{krasi.pdf}}%
  \end{picture}%
\endgroup%

%% file: Pictures/lifts.pdf_tex
\begingroup%
  \makeatletter%
  \providecommand\color[2][]{%
    \errmessage{(Inkscape) Color is used for the text in Inkscape, but the package 'color.sty' is not loaded}%
    \renewcommand\color[2][]{}%
  }%
  \providecommand\transparent[1]{%
    \errmessage{(Inkscape) Transparency is used (non-zero) for the text in Inkscape, but the package 'transparent.sty' is not loaded}%
    \renewcommand\transparent[1]{}%
  }%
  \providecommand\rotatebox[2]{#2}%
  \newcommand*\fsize{\dimexpr\f@size pt\relax}%
  \newcommand*\lineheight[1]{\fontsize{\fsize}{#1\fsize}\selectfont}%
  \ifx\svgwidth\undefined%
    \setlength{\unitlength}{517.65251656bp}%
    \ifx\svgscale\undefined%
      \relax%
    \else%
      \setlength{\unitlength}{\unitlength * \real{\svgscale}}%
    \fi%
  \else%
    \setlength{\unitlength}{\svgwidth}%
  \fi%
  \global\let\svgwidth\undefined%
  \global\let\svgscale\undefined%
  \makeatother%
  \begin{picture}(1,0.42788864)%
    \lineheight{1}%
    \setlength\tabcolsep{0pt}%
    \put(0,0){\includegraphics[width=\unitlength,page=1]{lifts.pdf}}%
    \put(0.85822017,0.29563719){\color[rgb]{0,0,0}\makebox(0,0)[lt]{\lineheight{10.25}\smash{\begin{tabular}[t]{l}$\alpha$\end{tabular}}}}%
    \put(0.13788433,0.26017813){\color[rgb]{0,0,0}\makebox(0,0)[lt]{\lineheight{10.25}\smash{\begin{tabular}[t]{l}$\widehat \alpha_2$\end{tabular}}}}%
    \put(0.04301155,0.2786717){\color[rgb]{0,0,0}\makebox(0,0)[lt]{\lineheight{10.25}\smash{\begin{tabular}[t]{l}$\widehat \alpha_1$\end{tabular}}}}%
    \put(0.73489471,0.38810672){\color[rgb]{0,0,0}\makebox(0,0)[lt]{\lineheight{10.25}\smash{\begin{tabular}[t]{l}$f(c_\infty)$\end{tabular}}}}%
    \put(0.14880809,0.36897849){\color[rgb]{0,0,0}\makebox(0,0)[lt]{\lineheight{10.25}\smash{\begin{tabular}[t]{l}$c_\infty$\end{tabular}}}}%
    \put(0,0){\includegraphics[width=\unitlength,page=2]{lifts.pdf}}%
    \put(0.83147046,0.18981381){\color[rgb]{0,0,0}\makebox(0,0)[lt]{\lineheight{10.25}\smash{\begin{tabular}[t]{l}$y$\end{tabular}}}}%
    \put(0.14485846,0.19233213){\color[rgb]{0,0,0}\makebox(0,0)[lt]{\lineheight{10.25}\smash{\begin{tabular}[t]{l}$\widehat y_2$\end{tabular}}}}%
    \put(0.00619151,0.18800126){\color[rgb]{0,0,0}\makebox(0,0)[lt]{\lineheight{10.25}\smash{\begin{tabular}[t]{l}$\widehat y_1$\end{tabular}}}}%
    \put(0,0){\includegraphics[width=\unitlength,page=3]{lifts.pdf}}%
    \put(0.48860137,0.38218937){\color[rgb]{0.52156863,0.52156863,0.52156863}\makebox(0,0)[lt]{\lineheight{10.25}\smash{\begin{tabular}[t]{l}$f$\end{tabular}}}}%
    \put(0.17712833,0.05701475){\color[rgb]{0.40392157,0.40392157,0.61568627}\makebox(0,0)[lt]{\lineheight{10.25}\smash{\begin{tabular}[t]{l}$\widehat J(f)$\end{tabular}}}}%
    \put(0.05711663,0.09219908){\color[rgb]{1,0.27843137,0}\makebox(0,0)[lt]{\lineheight{10.25}\smash{\begin{tabular}[t]{l}$K$\end{tabular}}}}%
    \put(0.63704439,0.09010599){\color[rgb]{1,0.27843137,0}\makebox(0,0)[lt]{\lineheight{10.25}\smash{\begin{tabular}[t]{l}$K$\end{tabular}}}}%
    \put(0.76124552,0.06126795){\color[rgb]{0.40392157,0.40392157,0.61568627}\makebox(0,0)[lt]{\lineheight{10.25}\smash{\begin{tabular}[t]{l}$\widehat J(f)$\end{tabular}}}}%
  \end{picture}%
\endgroup%

%% file: Pictures/rellenos.pdf_tex
\begingroup%
  \makeatletter%
  \providecommand\color[2][]{%
    \errmessage{(Inkscape) Color is used for the text in Inkscape, but the package 'color.sty' is not loaded}%
    \renewcommand\color[2][]{}%
  }%
  \providecommand\transparent[1]{%
    \errmessage{(Inkscape) Transparency is used (non-zero) for the text in Inkscape, but the package 'transparent.sty' is not loaded}%
    \renewcommand\transparent[1]{}%
  }%
  \providecommand\rotatebox[2]{#2}%
  \newcommand*\fsize{\dimexpr\f@size pt\relax}%
  \newcommand*\lineheight[1]{\fontsize{\fsize}{#1\fsize}\selectfont}%
  \ifx\svgwidth\undefined%
    \setlength{\unitlength}{224.85581646bp}%
    \ifx\svgscale\undefined%
      \relax%
    \else%
      \setlength{\unitlength}{\unitlength * \real{\svgscale}}%
    \fi%
  \else%
    \setlength{\unitlength}{\svgwidth}%
  \fi%
  \global\let\svgwidth\undefined%
  \global\let\svgscale\undefined%
  \makeatother%
  \begin{picture}(1,0.99858282)%
    \lineheight{1}%
    \setlength\tabcolsep{0pt}%
    \put(0,0){\includegraphics[width=\unitlength,page=1]{rellenos.pdf}}%
    \put(0.65836548,0.1929366){\color[rgb]{1,0.27843137,0}\makebox(0,0)[lt]{\lineheight{10.25}\smash{\begin{tabular}[t]{l}$A_i$\end{tabular}}}}%
    \put(0.20260352,0.14337024){\color[rgb]{0.42745098,0.42745098,0.62745098}\makebox(0,0)[lt]{\lineheight{10.25}\smash{\begin{tabular}[t]{l}$\widehat J(f)$\end{tabular}}}}%
    \put(0.38745954,0.85515833){\color[rgb]{0,0,0}\makebox(0,0)[lt]{\lineheight{10.25}\smash{\begin{tabular}[t]{l}$c_\infty$\end{tabular}}}}%
    \put(0.63576299,0.40892665){\color[rgb]{1,0.27843137,0}\makebox(0,0)[lt]{\lineheight{10.25}\smash{\begin{tabular}[t]{l}$\widehat A_i$\end{tabular}}}}%
    \put(0,0){\includegraphics[width=\unitlength,page=2]{rellenos.pdf}}%
  \end{picture}%
\endgroup%

%% file: main.bbl
\begin{thebibliography}{99}

\bibitem[CMR]{CMR-2004}
D. K. Childers, J. C. Mayer, J. T. Rogers (2004) \emph{Indecomposable continua and the Julia sets of polynomials, II}, Topology and its Applications, Elsevier

\bibitem[BM]{BM}
M. Bonk, D. Meyer (2017) \emph{Expanding Thurston maps}, ISBN-10: 0-8218-7554-X ISBN-13: 978-0-
8218-7554-4  

\bibitem[H]{hurwitz}
A. Hurwitz (1891) \emph{Uber Riemann'sche Fläche mit gegebenen Verzweigungspunkten.}, Math. Ann., 39 (1–60)

\bibitem[HY]{HY}
J. G. Hocking, G. S. Young (1961) \emph{Topology}, Addison-Wesley publishing company, London

\bibitem[IPRX]{IPRX 22}
J. Iglesias, A. Portela, A. Rovella, J. Xavier (2023)
\emph{Branched coverings of the sphere having a completely invariant continuum with infinitely many Wada lakes}, Topology and its Applications, 339, Part B

\bibitem[J]{janis}
Z. Janiszewski (1913) \emph{Sur les coupures du plan faites par les continus}, Prace Mat. Fiz. XXVI Warsaw

\bibitem[Kr72]{krasi 72}
J. Krasinkiewicz (1972) \emph{On the composants of indecomposable plane continua}, Bull. Acad. Polon. Sci. Sér. Math. Astronom. Phys., 20 (935 - 940)

\bibitem[Kr74]{krasi 74}
J. Krasinkiewicz (1974) \emph{On internal composants of indecomposable plane continua}, Fund. Math., 85 (255 - 263)

\bibitem[Ku]{kuratowski}
K. Kuratowski (1968) \emph{Topology vol. 2}, Academic Press and PWN, New York, London and Warszawa

\bibitem[MR]{MR-93}
J. C. Mayer, J. T. Rogers (1993) \emph{Indecomposable continua and the Julia sets of polynomials}, Proc. Amer. Math. Soc., Vol. 117, Num. 3 (795 - 802)

\bibitem[Mi]{milnor}
J. Milnor (2006) \emph{Dynamics in one complex variable, Third edition}, Annals of Mathematics
Studies, Princeton University Press

\bibitem[R]{rutt}
N. E. Rutt (1935) \emph{Prime ends and indecomposability}, Bulletin of the American Mathematical Society, Bull. Amer. Math. Soc., 41(4) (265 - 273)

\end{thebibliography}
